\input amstex
\documentstyle{amsppt}
\nologo
\parskip=6pt
\magnification=1200

\input epsf

\topmatter
\title Burau Representation and Random Walk on String
Links 
\endtitle
\author Xiao-Song Lin, Feng Tian and Zhenghan Wang
\endauthor

\abstract{Using a probabilistic interpretation of the Burau representation
of the braid group offered by Vaughan Jones, we generalize the Burau
representation to a representation of the semigroup of string links. 
This representation is determined by a linear system, and is dominated 
by finite type string link invariants.  For positive string links, the 
representation matrix can be interpreted as the transition matrix 
of a Markov process.  For positive non-separable links, we show that 
all states are persistent.
}
\endabstract

\thanks{Wang is supported by an NSF Postdoctoral Fellowship.}
\endthanks

\endtopmatter

\noindent{\bf \S1. Beginning of the story: Burau representation}
\medskip

The Burau matrices $\beta_i$, $1\leq i\leq n-1$, are $n\times n$ matrices 
given as follows: Let $t\neq 0,1$ be a complex number. If we think of
$\beta_i$ as a linear transformation on $\Bbb C^n$, then
$$\beta_i=\undersetbrace \text{$i-1$ copies} \to{(1)\oplus\cdots\oplus(1)}
\oplus\pmatrix 1-t & t \\ 1 & 0 \endpmatrix
\oplus\undersetbrace\text{$n-i-1$ copies}\to{(1)\oplus\cdots\oplus(1)},$$
and 
$$\beta_i^{-1}=\undersetbrace \text{$i-1$ copies} \to{(1)\oplus\cdots\oplus(1)}
\oplus\pmatrix 0 & 1 \\ \bar t & 1-\bar t \endpmatrix
\oplus\undersetbrace\text{$n-i-1$ copies}\to{(1)\oplus\cdots\oplus(1)}.$$
Here we use $\bar t$ to denote $t^{-1}$ for simplicity. 

It is easy to check that
$$\align& \beta_i\beta_{i+1}\beta_{i}=\beta_{i+1}\beta_{i}\beta_{i+1},\\
&\beta_i\beta_j=\beta_j\beta_i\qquad\text{for $|i-j|\geq2$}.
\endalign$$
Thus, sending the standard generators $\sigma_i$ of the braid group $B_n$ to 
$\beta_i$ defines the (non-reduced) Burau representation of $B_n$. 
For a given braid, its image under the Burau representation will be called
the {\it Burau matrix} of that braid.

There is an 
extensive literature on the Burau representation. We only mention an 
article of John Moody [7] where it was proved (settling a question of 
long time)
that the Burau representation is not faithful. In [5], Vaughan Jones 
offered a probabilistic interpretation of the
Burau representation. We quote from [3] (with a small correction):

``For positive braids there is also a mechanical interpretation of the Burau
matrix: Lay the braid out flat and make it into a bowling alley with $n$ lanes,
the lanes going over each other according to the braid. If a ball traveling 
along a lane has probability $1-t$ of falling off the top lane (and continuing
in the lane below) at every crossing then the $(i,j)$ entry of the 
(non-reduced) Burau matrix is the probability that a ball bowled in the $i$th 
lane will end up in the $j$th.'' 

Let us now consider string links. This notion was first introduced in [4].
We will generalize it a little bit here and still call the generalization 
string links. 

Essentially, a {\it string link} is an
{\it oriented tangle diagram} (or simply a {\it tangle}) 
in the strip $\Bbb R\times[0,1]$ with bottom ends $\{1\times 0,2\times0,
\dots,n\times0\}$ (call them {\it sources}) and top ends 
$\{1\times1,2\times1,\dots,n\times1\}$ (call them {\it sinks}).
 There are exactly $n$ strands,
each of them giving an oriented path from a bottom source $i\times 0$ to 
a top sink $j\times 1$. 
See Figure 1. Two such string links are thought to be the same if they 
differ by a finite sequence of
Reidermeister moves. Naturally, the set $S_n$ of all string links with 
$n$ strands has a semigroup structure such that $B_n\subset S_n$ is a 
subgroup.
 
\bigskip
\centerline{\epsfysize=1in \epsfbox{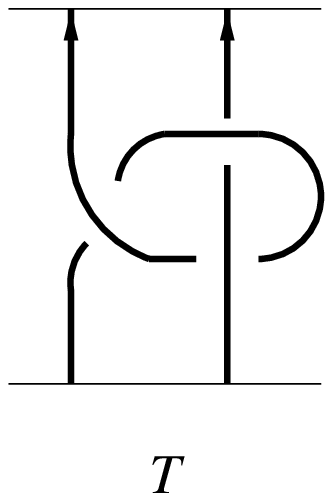}}
\smallskip
\centerline{$\ssize\text{Figure 1. A favorable string link.}$}
\medskip

We now define a representation of the semigroup $S_n$ generalizing
the Burau representation. We will assign to each element in $S_n$ 
an $n\times n$ matrix whose entries are rational functions in $t$.
Such an assignment will be multiplicative on $S_n$ so that we get a 
representation of the semigroup $S_n$ into the semigroup of $n\times n$ 
matrices.

Starting at the point $i\times0$,
we will try to walk up along strands of
the given string link $\sigma$ to get to the 
point $j\times1$ according to the the following rules:
\roster 
\item The walking direction should always 
be in agreement with the orientation of strands of $\sigma$. 
\item If we come to a crossing on the lower segment,
we will keep walking on the lower segment passing through that crossing.
\item If we come to a crossing on the upper segment,
we may choose either to jump down walking on the lower segment or 
keep walking on 
the upper segment passing through that crossing. 
\endroster 
Such a way of walking from $i\times0$ to $j\times 1$ on $\sigma$ is called a 
{\it path}.

A {\it loop} is a part of a path along which we may come back to where we
start on a string link. A path is called {\it simple} if it contains no 
loops. Obviously, 
there are only finitely many simple paths on a string link. A loop is {\it
simple} if it contains no other loops except itself. There are only finitely 
many simple loops on a string link. Any path can be reduced down to a simple
path by dropping off simple loops it contains. Therefore, there are at most 
countably many 
paths on a string link. 

Let us assign
a weight $w(P)$ to
a path $P$. Along $P$, there are many places where decisions are made 
about whether we jump down or keep walking. We will have a 
{\it state} at each of these places along $P$ and $w(P)$ is the 
product of all states on that path. The states are
determined as follows:
\roster
\item if we come to a positive crossing on the upper segment, 
the state is $1-t$ if we choose to jump 
down and
$t$ otherwise; and 
\item if we come to a negative crossing on the upper segment, the state
is $1-\bar t$ if we choose to jump down and $\bar t$ otherwise.
\endroster
 
With all these said, the $(i,j)$ entry of the $n\times n$ matrix assigned to 
the string link $\sigma\in S_n$ is   
$$ \sum_P w(P)\tag{1.1}$$
summing over all paths on $\sigma$ from
the point $i\times0$ to the point $j\times1$. It is zero if there is
no such a path.  We will call the matrix $B(\sigma)$ with $(i,j)$-th 
entry defined above the {\it Burau matrix} of $\sigma$. The invariance of $B(\sigma)$
under Reidermeister moves will be proved in Theorem B below.

This paper extensively revises the first author's preprint "Burau representation
and random walk on knots".  

\bigskip
\noindent{\bf \S2. Basic properties of the representation}
\medskip

In this section, we prove the following two basic theorems about 
the construction in \S1.

\proclaim{Theorem A} The sum {\rm (1.1)} converges to a rational function
of $t$.
\endproclaim

\proclaim{Theorem B} The matrix assigned to the string link $\sigma\in S_n$ is 
invariant
under Reidermeister moves.
\endproclaim

First we have the following observation.

We define the {\it multiplicity} of a path to be the number of simple loops
it contains. There are only finitely many paths with a given multiplicity. 

\proclaim{Observation:} 
If $P$ is a path with multiplicity $k$, then
$$w(P)=(1-t)^k\Phi(t)$$
where $\Phi(t)$ is a formal power series of $1-t$. 
The sum {\rm (1.1)} converges to a 
formal power series of $1-t$.  
\endproclaim

\demo{Proof} As formal power series of $1-t$, we have
$$\align
&\bar t=1+(1-t)+(1-t)^2+\cdots+(1-t)^n+\cdots\\
&1-\bar t= -(1-t)-(1-t)^2-\cdots-(1-t)^n-\cdots.
\endalign$$

Since a string link contains no closed components, the state 
of any simple loop contains a factor of $1-t$ or $1-\bar t$. This proves the 
first part. 

Now the second part follows easily from the fact that there
are only finitely many paths with a given multiplicity. \qed
\enddemo

Let $F^+=F^+_m$ be the free semigroup generated by $x_1,\dots,x_m$. Let
$\rho$ be a multiplicative function on $F^+$ such that $\rho(x_i)$ 
is a rational 
function of $t$ for every $i=1,\dots,m$. Moreover, assume that when $t=1$,
$\rho(x_i)=0$.

\proclaim{Proposition 2.1} The sum
$$\sum_{G\in F^+}\rho(G)$$
converges to a rational function of $t$. 
\endproclaim

\demo{Proof} We prove inductively a stronger statement: The sum
$$\Phi_m(t)=\sum_{G\in F^+_m}\rho(G)$$
converges to a rational function of $t$ such that $\Phi_m(1)=1$. Actually,
$$\Phi_m=\frac{1}{1-\dsize\sum_{i=1}^m\rho(x_i)}.\tag{2.1}$$
 
This statement is certainly true when $m=1$: Since $\rho(x_1)=0$ when $t=1$,
for $|1-t|$ small enough, $|\rho(x_1)|<1$ and therefore
$$\Phi_1(t)=\sum_{k=0}^\infty \rho(x_1^k)=\frac1{1-\rho(x_1)}.$$
Obviously, $\Phi_1(t)$ is a rational function of $t$ such that $\Phi_1(1)=1$.

We list elements in $F^+_m$ not equal to 1 as follows:
$$\align
&G_1x_m^{k_1}G_2x_m^{k_2}\cdots G_lx_m^{k_l},\\
&x_m^{k_1}G_1x_m^{k_2}G_2\cdots x_m^{k_l}G_l,\\
&x_m^{k_0}G_1x_m^{k_1}G_2x_m^{k_2}\cdots G_{l-1}x_m^{k_{l-1}},\\
&G_0x_m^{k_1}G_1x_m^{k_2}G_2\cdots x_m^{k_{l-1}}G_{l-1},
\endalign$$
where $l=1,2,\dots$, $k_j$'s are positive integers and $G_j$' are elements
in $F^+_{m-1}$ not equal to 1.

Inductively, we have
$$\align
\sum_{k_j,G_j,l}&\rho(G_1x_m^{k_1}G_2x_m^{k_2}\cdots G_lx_m^{k_l})\\
&=\sum_{l=1}^\infty\left[\left(\Phi_{m-1}-1\right)
\frac{\rho(x_m)}{1-\rho(x_m)}\right]^l.
\endalign$$
The last sum is convergent since $\Phi_{m-1}(1)-1=0$ and it is convergent to
a rational function of $t$. Moreover, the sum is zero when $t=1$. The same 
conclusions hold for other kinds of terms in the above list of elements
in $F^+_m$. Therefore, our inductive statement holds and this finishes the
proof of Proposition 2.1. Moreover, a straightforward calculation will 
establish $(2.1)$. \qed
\enddemo 

\demo{Proof of Theorem A} Multiple loops can be thought of as elements in
the free semigroup generated by simple loops. We already observed that 
the multiplicative weight function $w$ on each 
simple loop is zero when $t=1$. Since there are only finitely many simple
loops and finitely many simple paths, Theorem A follows immediately from 
Proposition 2.1. \qed
\enddemo

\demo{Proof of Theorem B}

There are three types of Reidermeister moves. We check the 
invariance one by one. Notice that in Figure 3, 4 and 5, the orientations of 
strands are determined by that of paths. 

\medskip
\noindent{\it Type I:} See Figure 2. For a positive kink, there are two paths 
(both simple) from the bottom end to the top end with 
$w=t$ and $w=1-t$ respectively. So 
total value is the same as no kink. For a negative kink, there is one simple
path from the bottom end to the top end with $w=\bar t$. There is also
a simple loop with $w=1-\bar t$. Therefore, the total value is
$$\sum_{k=0}^\infty\bar t(1-\bar t)^k=\frac{\bar t}{1-(1-\bar t)}=1.$$
Again, this is the same as no kink.

\medskip
\noindent{\it Type II:} See Figure 3. 
On the diagram on the left side of Figure 3, 
there are two paths going from the point $a$ to the point $b$ with
$w=1-\bar t$ and $w=\bar t(1-t)$ respectively. Their sum is
0, the same as no path from $a$ to $b$ on the diagram on the right side of 
Figure 3. Other situations can be checked similarly.

\medskip
\noindent{\it Type III:} See Figure 4. On the diagram on the left side of
Figure 4, there are two paths from $a$ to $b$ with $w=1-\bar t$ and 
$w=\bar t(1-t)^2=\bar t-2+t$ respectively. The total value is $t-1$, 
the same as the state
$t(1-\bar t)$ of the only path from $a$ to $b$ on the diagram on
the right side. Other situations can be checked similarly.

\bigskip
\centerline{\epsfxsize=1.7in \epsfbox{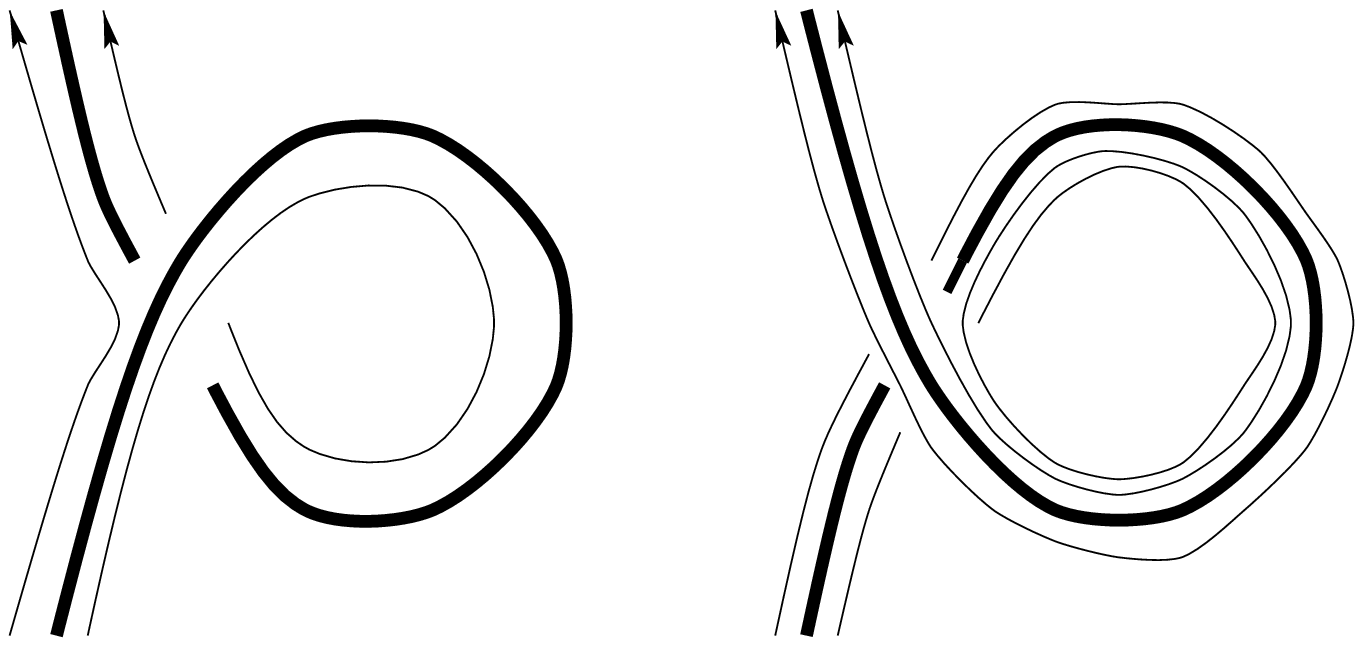}}
\smallskip
\centerline{$\ssize\text{Figure 2. Type I Reidemeister move.}$}
\medskip

\bigskip
\centerline{\epsfysize=1in \epsfbox{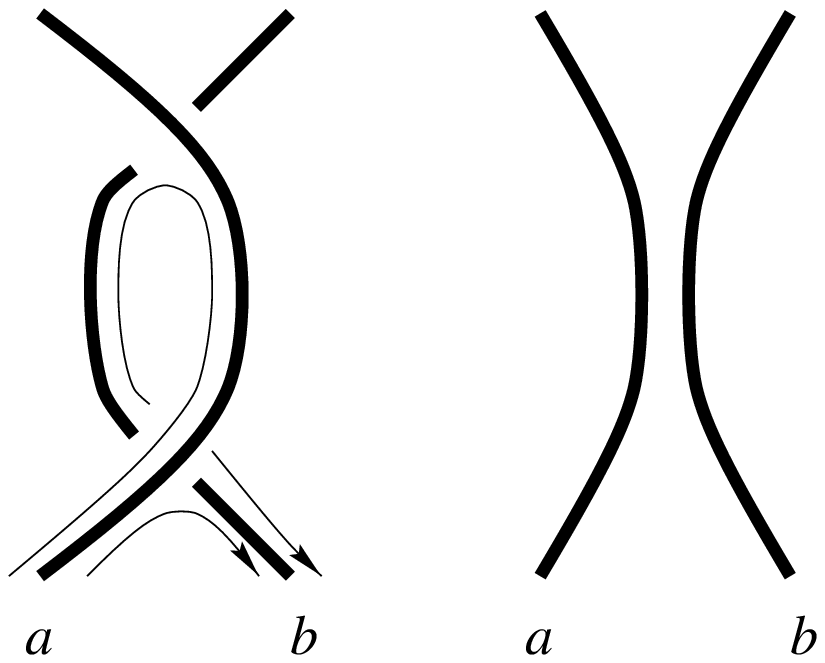}}
\smallskip
\centerline{$\ssize\text{Figure 3. Type II Reidemeister move.}$}
\medskip

\bigskip
\centerline{\epsfxsize=2.2in \epsfbox{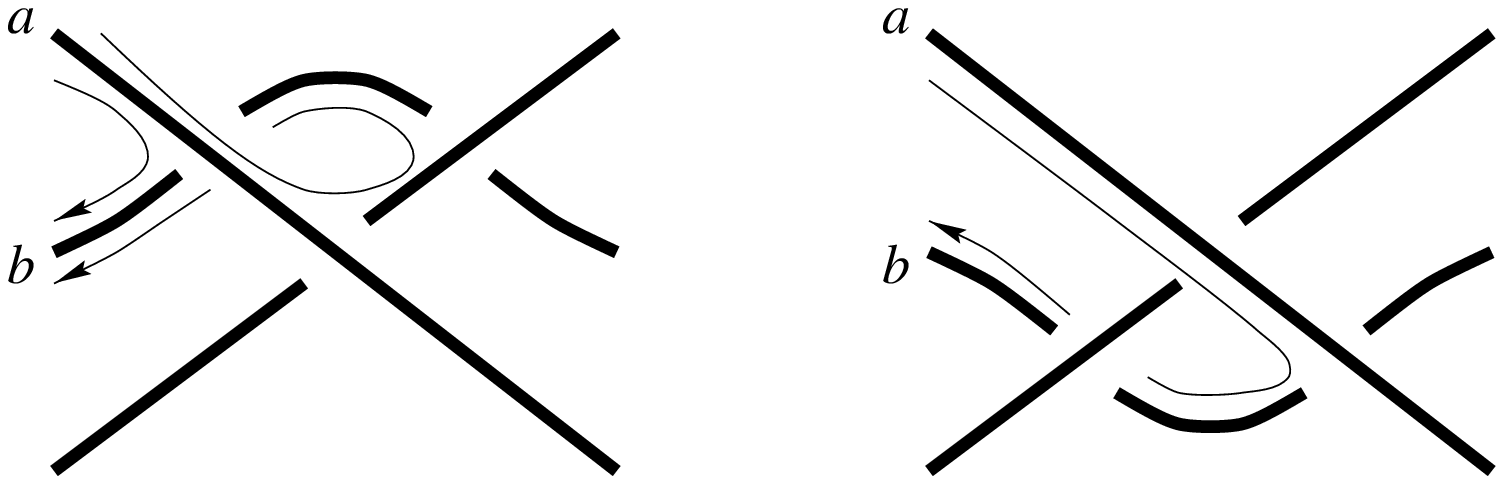}}
\smallskip
\centerline{$\ssize\text{Figure 4. Type III Reidemeister move.}$}
\medskip   
 
This finishes the proof of Theorem B.
\enddemo

\medskip
\demo{Example}
On the string link $T$ in Figure 1, there is only one simple loop with 
$w=t(1-\bar t)$. There is one simple path for the entry (1,1) with $w=1$,
one simple path for the entry (1,2) with $w=1-t$, and one simple path
for the entry (2,1) with $w=\bar t(1-t)$. All of them can be 
combined with multiple copies of the simple loop to form new paths. 
For the entry (2,2), there are two simple paths with $w=t$ and 
$w=(1-t)^2(1-\bar t)$ respectively. Only the second simple path
can be combined with multiple copies of the simple loop to form new paths.
The sum of the weights of all loops is convergent:
$$\sum_{k=0}^\infty[t(1-\bar t)]^k=\frac{1}{2-t}.$$
Therefore the Burau matrix of $T$ is
$$\frac{1}{2-t}\pmatrix 1 & 1-t \\ 
\bar t -1  & 3-t-\bar t \endpmatrix.$$
\enddemo

\bigskip
\noindent{\bf \S4. Path structures of string links}
\medskip

In this section, we discuss the path structure of string links. This 
leads to an algorithm to compute our representation.

Given a string link $\sigma \in S_{n}$, we may encode all paths on $\sigma$ from
a fixed bottom source to the top sinks into a binary tree. We first order
all overcrossings on the string link from the first component to the last one
following the orientations, and denote by $a_{kl}$ the $k$-th overcrossing on 
the $l$-th component.
Then all paths from the bottom $i$-th source to the top sinks are 
in one-one correspondence with paths from $a_{1i}$ in the following 
binary tree:
$$ \matrix  & & & \nearrow &  \cdots \\
												& & a_{2i} & & \\
             & \nearrow &  & \searrow & \cdots \\
          a_{1i} & & & &  \\
                 & \searrow &  &\nearrow & \cdots \\
            & & a_{kj} &  & \\
                 & & & \searrow & \cdots  \endmatrix\qquad,$$ 
where the $a_{ki}$ following an uphill edge is the next overcrossing  
when we keep walking, and the $a_{ki}$ following a downhill edge 
is the next overcrossing when we jump down. The trees end when we 
arrive at a sink after an uphill or downhill edge.  
We have a total of $n$ trees, and each path (directed) 
in the tree is finite, but there are infinitely many paths.

Suppose $x,y$ are two points on the string link.  Let $P(x,y)$ be 
the set of all paths from $x$ to $y$.  Since all activities happen 
at overcrossings, $P(x,y)$ and $P(x',y')$ are naturally isomorphic 
if there are no overcrossings between $x,x'$ and $y,y'$.  So 
in the following, we do not distinguish isomorphic path sets in this sense. In
the case when a point on the string link is denoted by $a_{kl}$, it is the point
on the $l$-th component right before the overcrossing $a_{kl}$.

Each path in $P(x,y)$ can be assigned a weight as before, and the total 
weight $wP(x,y)$ is the sum of the weights of all paths.   
Let us denote by $A_{kl}$ the row vector with $n$ components given 
by $(wP(a_{kl}, 1\times 1),\dots , wP(a_{kl}, n\times 1))$.  
We have the following theorem. 

\proclaim{Theorem C} The vectors $A_{kl}$ satisfy a linear system
in the field of rational functions of $t$, and 
this linear system has a unique solution around $t=1$.  
\endproclaim

\demo{Proof} The linear system is given by the binary tree as follows:
Before we arrive at the sinks, put $t^{\epsilon}$ on each uphill edge, and
$1-t^{\epsilon}$ on each  downhill edge,
where $\epsilon=\pm 1$ is the sign of the corresponding overcrossing.
At each vertex, say $a_{kl}$, followed by $a_{(k+1)l}$ (uphill), 
and $a_{mj}$ (downhill), 
the fundamental relation is 
$$ A_{kl}= t^{\epsilon} A_{(k+1)l} +(1-t^{\epsilon}) A_{mj}.$$
In this relation, if the vertex $a_{kl}$ goes to a sink following an uphill edge, 
$A_{(k+1)l}$ should be replaced by the vector $(0,\dots,0,1,0,\dots,0)$ where 1
appears at the
$l$-th position; and if $a_{kl}$ goes to a sink following a downhill edge,
$A_{mj}$ should be replaced by the vector $(0,\dots,0,1,0,\dots,0)$ where 1
appears at the $j$-th position. 

The existence of solutions follows from Theorem A.  To prove 
the uniqueness, note that when we solve the linear system by 
eliminations, we need only to multiply equations by 
$t^{\epsilon}$ or
$1-t^{\epsilon}$, and do substitutions. \qed 
\enddemo

As Burau matrix is 
given by $A_{1i}, i=1,\cdots ,n$ as rows, so the representation can be calculated 
by solving this linear system. 
For a string link with one component, all paths end at the same sink. 
So there are the obvious solutions $A_{kl}=1$.  Hence, we have 

\proclaim{Corollary 1:} If $\sigma \in S_{1}$, then $B(\sigma)=1$.
\endproclaim

In general, we get a new linear system by
summing all components in a vector.  This leads to the following 
important corollary. 

\proclaim{Corollary 2:} If $\sigma \in S_{n}$, then $B(\sigma)$
always has 1 as an eigenvalue with an eigenvector $(1,\cdots ,1)$.
\endproclaim

\noindent{\it Remarks:}

(1) If we exchange the role of sources and sinks, the representation matrix
will change dramatically. In this sense, we say that our representation depends on
the string link orientations. 

(2)  There are two obvious mirror images for string links: 
Imagine that a string link is laid almost flat in a horizontal plane. Then 
we have reflections in a vertical mirror or horizontal mirror, keeping the role
of sources and sinks unchanged. In the vertical case, the matrix 
will change by $t\mapsto\bar t$.  
The other case is complicated.

(3) Theorems A and B can also be proved
using Theorem C.  Actually, $t^{\epsilon}$ and 
$1-t^{\epsilon}$ are the only possible coefficients in the fundamental relations
which make the solution $A_{1i}$ invariant under Reidermeister moves.

We can design an electrical device using a 
positive string link (all overcrossings are positive):
we imagine that electrons enter from the bottom of a string link with a constant
density, and at each overcrossing the current is branched to 
upper and lower segments with $t$ and $1-t$ percents of the entering density,
respectively.  Then the sum of each column in the Burau matrix 
gives the output current density at the top.
If we reverse the process, so electrons 
enter from the top according to 
the output density function, we get a constant electric current. It seems to be an
interesting question whether or not the output current functions are dense in an
appropriate function space if we use all possible string links.  

\bigskip
\noindent{\bf \S5. Dominance by finite type invariants}
\medskip

It we are allowed to have some crossings of a string link replaced by transverse
double points, we get a singular string link. As usual, the Burau representation
matrix can be extended canonically to singular string links.  Given a string link 
$\sigma \in S_{n}$, and let $h=1-t$.  Then we can expand 
$B(\sigma)$ into $\sum_{i=0}^{\infty} b_{k}(\sigma) h^{k}$.  
Obviously, $b_{k}$ are also invariants of string links.
An invariant of string link is called of finite type of order $\leq k$ 
if it vanishes on any singular string link with more than $k$ double 
points. See [1,2,6].

\proclaim{Theorem D:}  The invariants $b_{k}$ is of finite type 
of order $\leq k$.
\endproclaim 

\demo{Proof} Let $\tau$ be a singular string link with $n>k$ double points
and let $\sigma$ be the string link obtained from $\tau$ by resolving all double
points to overcrossings. A path $P$ on $\sigma$ will induce paths on other
resolutions of $\tau$, having zero weight if the induced path is illegal. The
alternating sum of weights of all these induced paths is called the contribution of
$P$ in
$B(\tau)$. Thus, if a path on
$\sigma$ does not pass through some double points on $\tau$, its contribution in
$B(\tau)$ is zero. Otherwise, its contribution is zero modulo $h^n$. This could be
seen easily as the weight difference of an overcrossing and an undercrossing is
divisible by $h$.
\qed
\enddemo

This dominance of $B(\sigma)$ by finite type invariants is
closely related with the following interesting 
fact: Although there might be infinitely many paths from 
a source to a sink, there are only finitely many with 
nonzero weights modulo $h^{n}$.  In other words, 
if we restrict the number of jumping downs, there are only 
finitely many paths from a source to a sink. 

We don't know how to related our representation with the universal representation
coming from iterated integrals [6]. 

\bigskip
\noindent{\bf \S6. Markov processes on positive string links}
\medskip

Positive string links have only positive crossings. For such string links, 
the Burau
representation could be used to determine Markov
processes. 

First let us recall some basic notions about Markov processes.
All facts can be found in [3].

Let $S$ be a finite or countable set.
Suppose that to each pair $i,j$ in $S$ that is assigned a nonnegative 
number $p_{ij}$, and that these numbers satisfy the equations:
$$ \sum_{j\in S} p_{ij}=1, \;\; i\in S. \tag{6.1}$$
Let $X_{0}, X_{1}, \dots $ be a sequence of random variables which take values in
$S$.  The sequence is a {\it Markov process} if the conditional
distribution  of the next state $X_{n+1}$ given the present state $X_{n}$ is
independent  of the past states $X_{0},\cdots , X_{n-1}$.  The elements of $S$ are 
thought as the possible states of a system, $X_{n}$ representing the state 
at time $n$.  The sequence or process $X_{0}, X_{1}, \dots$
then represents the history of the system, which evolves in accordance with 
a probability law determined by the transition matrix $P=(p_{ij})$, and 
the initial probabilities $\pi_{i}=P[X_{0}=i]$.

A matrix whose entries are nonnegative and satisfy equation (6.1) is 
called a {\it stochastic matrix}.  By an existence theorem (see [3], Theorem
8.1), given a stochastic matrix $P=p_{ij}$ and nonnegative numbers $\pi_{i}$
satisfying $\sum_{i}\pi_{i}=1$, there exists a Markov process on some 
probability space with initial probabilities $\pi_{i}$ and transition 
matrix $P$.

Given a positive string link, Corollary 2 in {\bf \S 4} tells us that 
for each $t\in (0,1]$, the representation matrix is a stochastic matrix.
Hence, given any initial conditions, there exists a Markov process with 
the components of the string link as the state space, and the transition matrix 
$B(\sigma)$.  So we may imagine that a particle travels on a string link, 
and the state is the component it stays on.  Each step it goes 
from one component (from bottom) to another (at top) with the 
probability equal to the total weight of all possible paths. 

Let $p_{ij}^{(n)}$ be the $(i,j)$-th entry of the matrix $P^{n}$.
A state $i\in S$ is {\it persistent} if $$\sum_{n} p_{ii}^{(n)} =+\infty.$$
Otherwise it is {\it transient}. 
All states are persistent if $$\sum_{n} p_{ij}^{(n)}=+\infty$$ for all $i,j$.
In this case, the system visits every state infinitely often.

A string link is {\it separable} if 
the link obtained by connecting the  
$i$-th source to the $i$-th sink (i.e., the usual closure of string links) 
can be separated by a 2-sphere into two non-empty links.

\proclaim{Theorem E}  Given a non-separable positive string link, 
then for each $t\in (0,1)$ 
all states of the associated Markov process are persistent. 
\endproclaim

We first have the following 
two lemmas.

\proclaim{Lemma 1:}  If $ \sigma \in S_{n}$ is a non-separable  
string link, then there exists an integer $N$ such that  for every pair 
$i,j$, there is a path from the $i$-th source to the $j$-th sink in 
${\sigma}^{N}$. 
\endproclaim

\proclaim{Lemma 2:} Suppose $A$ is a stochastic matrix.  If there exists 
an integer $N$ such that $a_{ij}^{(N)}>0$ for every pair $i,j$. Then 
\roster
\item 1 is the unique non-negative eigenvalue of $A$ with a strictly 
positive (row or column) eigenvector.
\item No eigenvalue of $A$ has absolute value greater that 1.
\item Let $u=(u_{1}, \dots , u_{n})$ be the unique row (left) eigenvector of $A$,  
then for each pair $i,j$,  $$\lim_{n\rightarrow\infty} a_{ij}^{(n)}=u_{j}.$$
\endroster
\endproclaim

Lemma 2 follows from Perron-Frobenius theory and the renewal theorem (see pp. 16--17
of [8]). We first prove Theorem E, assuming Lemma 1 and Lemma 2.

\demo{Proof of Theorem E} By Lemma 1, 
for each pair $i,j$, we have a path from the $i$-th source to the 
$j$-th sink on $\sigma^N$.  So for $t\in (0,1)$, the $(i,j)$-th entry in
$B(\sigma^{N})$ is positive.  
So, by Lemma 2, 
$$\lim_{n\rightarrow\infty} p_{ij}^{(n)}=u_{j}>0$$
for each pair $i$ and $j$.
It follows that $$\sum_n p_{ij}^{(n)}=+\infty.$$ 
This finishes the proof of Theorem E. \qed
\enddemo

\demo{Proof of Lemma 1} We proceed by induction.  The case $n=1$ is trivial.  

Suppose our Lemma holds for $\leq n$ strands. For a string link $\sigma$ with
$n+1$ strands, taking some power if necessary, we may assume that the
permutation on $n+1$ letters induced by $\sigma$ is trivial so that the $i$-th
strand connects the $i$-th source to the $i$-th sink. (Note: We don't claim 
that any power of a non-separable string link is still non-separable, although
this seems to be true.)  
Delete the first strand and
we get a string link with $n$ strands. The strands of this string link 
are partitioned into non-separable string  
links $\Delta_k$ with $\leq n$ strands. Since the induced permutation is 
the identity, powers of $\sigma$ will be partitioned into powers of
$\Delta_k$. Now by taking some power of the string link $\sigma$ if
necessary, we may assume inductively that each string link $\Delta_k$ has the
property that any source and any sink is connected by a path. Since the
string link $\sigma$ is non-separable, 
the first strand has to have overcrossings
and undercrossings with strands in each $\Delta_k$. 
Then, it is easy to argue that
on $\sigma^N$ for some $N$, 
any source and any sink can be connected by a path. This finishes
the induction. \qed

\enddemo

\noindent{\it Remarks:}

(1)  As any link decomposes into non-separable links, it follows 
that for any positive string link, 
each state $i$ is persistent.

(2)  When $t=1$, the stochastic matrix is a permutation matrix.
When $t=0$, the limiting matrix is stochastic.  But it seems to be hard to
make a general statement.

(3)  There is no difficulty to consider infinite string links.  The situation 
is much like a 1-dimensional random walk.  But they differs essentially.
First our associated Markov process does not have spatial homogeneity.
Secondly, the transition matrix is not local.  It is possible for 
the particle to walk into any other position.  It makes more sense to 
think of it as a "quantum 1-dimensional random walk".

(4)  The existence of a row (left) strictly positive eigenvector 
corresponding to eigenvalue 1 makes the Markov process into a stationary 
distribution.  So we can then define the entropy for the string link for each 
$t$ as in information theory.  It would be interesting to relate 
entropy with the complexity of the string link.

\bigskip

\noindent Xiao-Song Lin
\newline The Chinese University of Hong Kong and 
\newline University of California, Riverside
\newline E-mail: {\it xslin\@ims.cuhk.edu.hk} or {\it xl\@math.ucr.edu}
\medskip

\noindent Feng Tian
\newline University of California, San Diego 
\newline E-mail: {\it ftian\@euclid.ucsd.edu}
\medskip
 
\noindent Zhenghan Wang 
\newline University of California, San Diego and 
\newline Indiana University
\newline E-mail: {\it zwang\@euclid.ucsd.edu}
\enddocument